\documentclass[a4paper,12pt]{amsart}

\usepackage{amssymb,amsbsy,amsmath,amsfonts,amssymb,amscd}
\usepackage{latexsym}
\usepackage{graphics}
\usepackage{color}
\input xy
\xyoption{all}

\theoremstyle{plain}
\newtheorem{thm}{Theorem}[section]

\newtheorem{lemma}[thm]{Lemma}

\theoremstyle{definition}

\newtheorem{defin}[thm]{Definition}

\newtheorem{example}[thm]{Example}

\newtheorem{question}[thm]{Question}

\numberwithin{equation}{section}

\newcommand{\sA}{{\mathcal A}}
\newcommand{\sB}{{\mathcal B}}
\newcommand{\sC}{{\mathcal C}}

\newcommand{\sF}{{\mathcal F}}

\newcommand{\sM}{{\mathcal M}}

\newcommand{\sU}{{\mathcal U}}

\newcommand{\PP}{\ensuremath{\mathbb{P}}}

\newcommand{\CC}{\ensuremath{\mathbb{C}}}
\newcommand{\RR}{\ensuremath{\mathbb{R}}}
\newcommand{\ZZ}{\ensuremath{\mathbb{Z}}}
\newcommand{\QQ}{\ensuremath{\mathbb{Q}}}

\newcommand{\NN}{\ensuremath{\mathbb{N}}}
\newcommand{\hol}{\ensuremath{\mathcal{O}}}

\newcommand\om{\omega}
\newcommand\la{\lambda}
\newcommand\s{\sigma}

\newcommand\ga{\gamma}

\newcommand\e{\epsilon}

\newcommand{\ra}{\ensuremath{\rightarrow}}

\def\eea{\end{eqnarray*}}
\def\bea{\begin{eqnarray*}}

\newcommand\dual{\mathrel{\raise3pt\hbox{$\underline{\mathrm{\thinspace d
\thinspace}}$}}}
\newcommand\qe{\ifhmode\unskip\nobreak\fi\quad $\Box$}       

\def\BOX{\hfill\lower.5\baselineskip\hbox{$\Box$}}

\newtheorem{theo}{Theorem}[section]

\newtheorem{remarkk}[theo]{Remark}
\newenvironment{rem}{\begin{remarkk}\rm}{\end{remarkk}}

\newtheorem{prop}[theo] {Proposition}
\newtheorem{cor}[theo]{Corollary}

\newenvironment{ex}{\begin{example}\rm}{\end{example}}

\newcommand{\Proof}{{\it Proof. }}

\title [Canonical maps of SIP's]{On the canonical map of some surfaces isogenous to a product}

\author{Fabrizio Catanese}
\address {Lehrstuhl Mathematik VIII\\
Mathematisches Institut der Universit\"at Bayreuth\\
NW II,  Universit\"atsstr. 30\\
95447 Bayreuth}
\email{fabrizio.catanese@uni-bayreuth.de}

\thanks{
The present work took place in the framework  of the 
 ERC Advanced grant n. 340258, `TADMICAMT' }

\date{\today}

\begin{document}

\maketitle

{\em  Dedicated to Lawrence Ein  on the occasion of his 60th birthday}

\begin{abstract}
We construct several families (indeed, connected components of the moduli space) of surfaces $S$ of general type with $p_g=5,6$ whose canonical map has image $\Sigma$ of very high degree,
$d=48$ for $p_g=5$, $d=56$ for $p_g=6$. And a connected component of the moduli space consisting of surfaces $S$ with $K^2_S = 40, p_g=4, q=0$
whose  canonical map has always degree $\geq 2$, and, for the general surface,  of degree $2$ onto
a canonical surface $Y$ with $K^2_Y = 12, p_g=4, q=0$.
 
The surfaces we consider are SIP 's, i.e. surfaces $S$ isogenous to a product of curves $(C_1 \times C_2 )/ G$;   in our examples the group $G$
is  elementary abelian, $G \cong (\ZZ/m)^k$.
We also establish some basic results concerning  the canonical maps of any surface isogenous to a product, basing on elementary representation theory.

Our methods and results are a first step towards   answering  the  question of existence of SIP 's  $S$ with $p_g=6, q=0$ whose canonical map  embeds $S$
as a surface of degree $56$ in $\PP^5$. 
\end{abstract}

\tableofcontents

\section{Introduction and state of the art on canonical surfaces of high degree in $\PP^3, \PP^4, \PP^5$}

Let $S$ be a complex minimal projective surface of general type, and denote by $p_g (S) : = dim H^0(\hol_S(K_S))$,
and by $\Phi $ the rational map to $\PP^{p_g-1}$ associating to a point $x \in S$ the point in $\PP (H^0(\hol_S(K_S)))^{\vee}$
obtained by evaluation of sections of $H^0(\hol_S(K_S))$  in $x$.

If the image of $\Phi$ is a surface $\Sigma$ (this is  the case if  the canonical system $|K_S|$ has no base points) 
one has a well known inequality:
$$ K^2_S \geq deg(\Phi) deg (\Sigma)  ,$$ 
equality holding iff $|K_S|$ has no base points.

\begin{defin}
We define the {\bf canonical degree} of $S$ as  $d : = deg (\Sigma)$ if $\Sigma$ is a surface, $d=0$ otherwise, and we define  the {\bf canonical volume} of $S$ to be $ K^2_S $.
Hence $ K^2_S \geq d$, equality holding if and only if $|K_S|$ has no base points and $\Phi $ is birational onto its image
(i.e., $deg (\Phi)=1$).

$S$ is said to be a {\bf (simple) canonical surface } if $\Phi $ is birational onto its image $\Sigma$. This obviously  implies that $p_g(S) \geq 4$.
\end{defin}

\begin{question}\label{q1}
What is the maximal canonical degree $d$  for a fixed value of $p_g$? 

What is the maximal canonical volume $K^2_S$?

In particular, for $p_g = 4,5,6$?
\end{question}

Recall the Castelnuovo inequality, holding if $\Phi$ is birational (onto its image $\Sigma$):
$$  K^2_S \geq 3 p_g(S) - 7,$$ 
and the  Bogomolov-Miyaoka-Yau inequality
$$  K^2_S \leq 9 \chi(S) = 9 + 9 p_g(S) - 9 q(S),$$
where $ q : = dim H^1 (\hol_S) $ (if $q=0$ $S$ is said to be {\bf regular}).

By virtue of these inequalities,  under the assumptions of  question \ref{q1} one must have:
 $$ d  \leq K^2_S \leq 9 (1 + p_g).$$
 In particular, the upper bound for the volume is $K^2_S \leq 45, 54, 63$ for $p_g = 4,5,6$.
 
 In order to have high volume it is better to have $q(S) = 0$, even if  also the canonical maps of irregular surfaces are a very interesting subject 
 of investigation, see \cite{c-s}. For instance Cesarano in his Ph.D. work considers surfaces $S$ which are a polarization of type $(1,2,2)$
 on an Abelian threefold, which are canonical surfaces in $\PP^5$ with $K^2_S= 24$, and he seems to have shown that for a general such surface
 the canonical map is an embedding. 
 
 \begin{question}\label{q2}
What is the maximal canonical degree $d$  of a surface for which $\Phi$ is an embedding?

In particular, for $p_g = 6$?

 (Note that in case of an embedding $ d = K^2$).

\end{question}

Indeed,  question \ref{q2} is easily answered for $p_g=4$, here only the case $d=5$
is possible, while,  for  $p_g=5$, must be $d=8,9$, and in all the three cases $\Sigma$ is a complete intersection,
of respective types $(5), (4,2), (3,3)$,
see \cite{AMS}, propositions 6.1 and 6.2, corollary 6.3., and \cite{Ellia}.
 
 Concerning the case $p_g=6$, in  the article \cite{AMS} methods of homological algebra were used to construct embedded canonical  surfaces with low degree $11 \leq K^2_S \leq 17$,
 and also   to attempt a classification of them. Recently, M. and G. Kapustka constructed in \cite{kk} such canonical surfaces of degree $K^2_S =18$,
 using the method of bilinkage.
 In \cite{Ellia} we constructed a family of  embedded canonical surfaces with $d=24, q=0$.
 
 Concerning the case $p_g=4$, against Enriques' expectation $K^2_S \leq 24$,   in \cite{BC} we constructed a
 surface with $p_g(S) = 4, q(S) = 0, K^2_S = 45$ (and with $K_S$ ample being a ball quotient): but here the canonical system has 
 base curves, $S$ is canonical and we have  $d=19$.
Surfaces with $p_g(S) = 6, q(S) = 2, K^2_S = 45 > d$ (again ball quotients) were constructed in \cite{BC}, but the canonical degree $d$ was not
calculated.

For $p_g=4$, the current record for the canonical degree is $d=28$ \cite{singular}, and all possible degrees in the interval
$ [5, 28]$ occur: Ciliberto got the interval $[5,10]$ in \cite{ciro}, Burniat the interval $[11,16]$ in \cite{burniat},
while the interval $[11,28]$ was gotten in \cite{singular}. Liedtke \cite{Liedtke} got $K^2_S = 31$, but $d=12$.

In this paper we construct surfaces with very high canonical degrees, $p_g=5, d = 48$, and  $p_g=6, d = 56$.

Indeed, the maximum degree a priori possible $d = 9 (p_g +1)$  can only occur for ball quotients with $q=0$, but these are difficult to
construct, and even more difficult it is to describe their canonical map.

We obtain $d = 8 (p_g +1)$ using  surfaces isogenous to a product with $q=0$,  for these the volume is $K^2 = 8 (p_g +1)$,
and the only difficulty is to describe their canonical maps.

Recall that surfaces isogenous to a product are quotients $(C_1 \times C_2)/G$ for the free action of a finite group $G$.
Their classification is a hard enterprise (but it was achieved for $p_g=0$ in \cite{bcg}, see also \cite{frapporti}).

In our paper we consider  easier cases where the group $G$ is Abelian (but we do not attempt a complete classification of the case $G$ Abelian and $p_g = 4,5,6$).
We lay however the foundations for the description of the canonical maps of surfaces isogenous to a product.

For this purpose a first step that we achieve is to decide when does a subrepresentation $ V \subset \CC[G]$ (i.e., a $\CC[G]$-submodule)
 yield a projective embedding
of the finite group $G$.
Another  open question is:
\begin{question}\label{q4}
Does there exist a surface $S$ isogenous to a product and with $p_g(S) = 6$
such that the canonical map $\Phi : S \ra \PP^5$ is an embedding?
\end{question}

Our main results in this article are the following ones, here written without full details (these will be given in the later sections).

\bigskip
 
 {\bf Theorem \ref{birat}}
 {\em 
 
I) There exists a connected component of the moduli space of the minimal surfaces 
of general type with $p_g=6, q=0, K^2_S = 56$, consisting of surfaces isogenous to a product with 
 $ G = (\ZZ/2)^3, $ 
such that the canonical map for each surface $S$
is a non injective birational morphism  onto a  canonical image $\Sigma$ of degree $56$ in $\PP^5$.

II) There exists a connected component of the moduli space of the minimal surfaces 
of general type with $p_g=5, q=0, K^2_S = 48$, consisting of surfaces isogenous to a product with 
 $ G = (\ZZ/2)^3, $ 
such that 
 the canonical map of $S$
is a non injective birational morphism  onto a  canonical image $\Sigma$ of degree $48$ in $\PP^4$,
for general choice of the branch points for $C_2 \ra \PP^1$.

}
\bigskip

There are other families of surfaces, beyond those mentioned in theorem \ref{birat}, that achieve the same $d$. The
families are here indicated, but the details are not worked out for them.

\bigskip

The next theorem provides  new features of an example which in the past had  not been considered possible (see \cite{Babbage}, \cite{beacan}, \cite{babbage}):
a canonical map of higher degree onto a canonical surface.  Previous examples (sometimes the same) appeared in \cite{vGZ}, \cite{beacan}, \cite{babbage},
\cite{pardini}, \cite{tan},  \cite{schoen}, \cite{tale}. In the previous examples the canonical map $\Phi$ was ramified only on a finite set,
and $\Phi$ was birational for a general deformation of
 the surfaces.

{\bf Theorem \ref{doublecan}}
{\em 
There is a connected component of the moduli space of the minimal surfaces 
of general type with $p_g=4, q=0, K^2_S = 40$, consisting of surfaces isogenous to a product with 
$G = (\ZZ/2)^3$,  such that for each of them there is a nontrivial involution $\psi : = (e,1) \in G \times G$ which acts  trivially on the canonical
system, and with quotient a nodal  surface $Z$ with $K_Z^2= 12$, and with $K_Z$ ample.

The canonical map of $S$ factors through the canonical map of $Z$,
which  is a morphism, and is birational onto its
image for general choice of the six branch points of $C_1 \ra  C_1 / G = \PP^1$.

The fixed locus of $\psi$ consists of $14$ isolated points, and two curves $F'_1 + F'_2$ which form the  base locus  of the canonical system .}

 \bigskip
 
 The last  result we mention in the introduction is meant   to give a flavour of the methods developed  in section $2$:

\smallskip 
 {\bf Proposition \ref{projemb}}
 {\em   A subrepresentation $V \subset \CC[G]$ yields a map of  $G$ into $\PP (V^{\vee})$  iff $V \neq 0$; it embeds $G$  only if the following properties
   i) and ii) hold.

  i)  Writing $  V =  \sum_{\chi \in Irr}  V_{\chi }, V_{\chi } = W_{\chi } \otimes U_{\chi},$
the representation $$\rho : = \otimes_{ U_{\chi}\neq 0} \  \rho_{\chi}$$
   is faithful.
   
   ii) Let $\chi_1, \dots , \chi_m$ be the irreducible representations for which $ U_{\chi} \neq 0$, and notice  that  $\rho_j$ sends  the centre $Z(G)$ of $G$
to the set of scalar  multiples of the identity,  via a character $\psi_j$ of $ Z(G)$.  Then the characters $\psi_1, \dots, \psi_j$ are affine generators of  
the group of characters $Z(G)^*$
(i.e.,  $\frac{\psi_j }{\psi_1}$ generate the group $Z(G)^*$).

If moreover $V$ is full (each space  $ U_{\chi} \neq 0$  is $= W_{\chi }^{\vee}$), the two conditions i) and ii) are also sufficient.

If $V$ is not full, we obtain a projective embedding of $G$ if $V$ is projectively general (i.e., each space  $ U_{\chi} \neq 0$  is projectively general).
  }

\section{ Representations yielding an affine, respectively   a  projective embedding, of a finite group $G$}

In this section $G$ is a finite group, and $V$ is a finite dimensional complex vector space which is a  representation of $G$ (indeed, one could  
replace $\CC$ by  an algebraically closed  field $K$ of characteristic prime to the order of $G$).

Assume that we have a linear map of representations
$$ \e : V \ra \CC[G] = \{ \sum_{g \in G} a_g g| a_g \in \CC \}, $$
where  we view the group algebra $\CC[G]$ as the space of complex valued functions on $G$, and $g$ denotes for brevity  the characteristic function of an element $g \in G$.
 Equivalently, $V$ is a left  $\CC[G]$-module and $\e$ is a homomorphism of $\CC[G]$-modules.

\subsection{Embedding $G$}
The first question is: when does $V$ separate points of $G$? I.e., we ask when  is the map $\e^* : G \ra V^{\vee}$,  obtained composing
 the inclusion $G \subset  \CC[G]^{\vee} $ with the dual map $\e^{\vee}$, injective?

The second question is: when does $\PP(\e^*)$   yield a projective embedding of $G$ inside $\PP(V^{\vee})$?

Replacing $V$ by its image inside $ \CC[G]$, we may assume in this section,  without loss of generality, that $ V \subset \CC[G]$.

A useful  observation is that if $ V \subset V' \subset \CC[G]$ and $V$ separates points of $G$, then so does a fortiori $V'$;
and similarly for the projective  embedding since $V'^{\vee} \ra V^{\vee}$ is surjective.

We denote by $\sA $ the group algebra $\CC[G]$; its centre $Z_{\sA} : = Z (\CC[G])$ is the space of class functions, the complex vector space 
with basis $\{ v_{\sC}\}$, indexed by the conjugacy classes $\sC$ of $G$, and where
$$ v_{\sC} : = \sum_{g \in \sC} g.$$

Another more useful basis is  indexed by the irreducible complex representations $W_{\chi}$ of $G$,
and their characters $\chi$ (the irreducible characters  form an orthonormal basis for the space of class functions).
We shall, by the way, denote by $\rho_{\chi} : G \ra GL (W_{\chi})$ the corresponding representation homomorphism,
and observe that there is a canonical isomorphism (dual map) $ GL (W_{\chi}) \cong  GL (W_{\chi}^{\vee})$.

One associates as usual, to an irreducible character $\chi$, an idempotent in the group algebra:
$$ e (\chi) : = \frac{\chi(1)}{|G|} \sum_{ g \in G} \chi (g^{-1}) g.$$ 

The elements $ e (\chi) $ are the principal  idempotents of the group algebra $\CC[G]$, and,  from the fact that  their sum  is
equal to $1$, follows that the group algebra is semisimple,
indeed

$$\sA = \CC[G] = \sum_{\chi \in Irr} \CC[G] e (\chi) = :   \sum_{\chi \in Irr} \sA{\chi}.$$

By the same token,  each  representation $V$ splits as a direct sum of irreducible representations,
$$  V =  \sum_{\chi \in Irr}  V_{\chi }, V_{\chi } = W_{\chi } \otimes (\CC^{n(\chi)}),$$
where $ \CC^{n(\chi)}$ is a trivial representation of $G$
(the orthogonality of characters imply that
$ e (\chi) $ acts as the identity on $V_{\chi }$ and acts as zero on $V_{\chi' }$ for $\chi' \neq \chi$).

 The lemma of Schur says that   $Hom_{\sA} (W_{\chi }, W_{\chi' })=0$ for $\chi' \neq \chi$,
and $ \cong \CC$ for $\chi'  =  \chi$.

   If $V$ is a left  $ \CC[G]$-module,  $V^{\vee} :=  Hom_{\CC} (V, \CC)$
is a right $ \CC[G]$-module for the (duality preserving) operation such that 
$$ v (f g ) : = (g^{-1} v) f, \  v \in V, \ f \in V^{\vee} \Rightarrow (gv) (fg) = v f.$$ 
Since there exists on $V$ a $G$-invariant positive definite Hermitian product,
the matrices of elements of $G$ in a unitary basis are unitary,
so that, viewing $V^{\vee}$ as a left $ \CC[G]$-module for the operation
$ f \mapsto f g^{-1}$, the character of $V^{\vee}$ is the conjugate of the one of $V$.

Moreover,  $$ \CC[G] = \sum_{\chi \in Irr} \sA{\chi}, \  \sA{\chi} \cong End (W_{\chi})= W_{\chi} \otimes W_{\chi}^{\vee}.$$
And, in the last isomorphism, $ e (\chi) $ corresponds to the identity tensor.
 In fact, $ \sA = \CC[G]$ is also a right $ \CC[G]$-module, and 
$ \sA _{\chi } = e (\chi ) \sA  = \sA e (\chi)$, so $ \sA _{\chi } $ is naturally of the form
$W_{\chi} \otimes W'$ and of the form $W''  \otimes W_{\chi}^{\vee}$.
Since the  dimension of $ \sA _{\chi }$ equals the square of the dimension of $W_{\chi} $,
then necessarily $ \sA _{\chi } = W_{\chi} \otimes W_{\chi}^{\vee}.$

\begin{rem} Replacing $\CC$ by any other field $K$ of characteristic zero, one sees that the semisimple decomposition
of $K[G] = : \sA_{K}$ is given by the direct sum over  $\sA_{K} e_K(\chi)$,
where $$ e_K(\chi): = \sum_{\sigma \in Gal (K(\chi), K)} e (\chi^{\sigma})$$
 and where $K(\chi)$ is the extension of $K$ generated
by the values of the character $\chi$.

Again,  $\sA_{K} e_K(\chi)$ is a matrix algebra over a division algebra $D_{\chi}$.

\end{rem}

\begin{defin}\label{full}
1) We shall say that a left $\CC[G]$- submodule $ V \subset \CC[G]$ is a  {\bf full } subrepresentation if it is a direct sum $\sA{\chi_1}\oplus \dots \oplus \sA{\chi_m}$.

This holds automatically if $G$ is Abelian, since in this case each irreducible representation is $1$-dimensional.

2) In general, for $ V \subset \CC[G]$,  $V_{\chi } = W_{\chi } \otimes U_{\chi }$, for $U_{\chi }$ a vector subspace of $W_{\chi}^{\vee}.$

2.1) We shall say that $U_{\chi }$ is {\bf  linearly general } if,   setting  $G' : = \rho_{\chi} (G)$,  $U_{\chi}$ is not contained in any  of the proper subspaces
   $$  Ker  (\ga - \ga'),  \ for \  \ga \neq \ga' \in G''.$$ 
  (then it is not contained in their union).
  
  We shall say that $V$  is {\bf  linearly general } if every $U_{\chi } \neq 0$ is   linearly general.

2.2) We shall say that $U_{\chi }$ is {\bf  projectively general } if,   setting  $$G'' : = \PP (\rho_{\chi} (G)) \subset \PP GL (W_{\chi } )= 
\PP GL (W_{\chi }^{\vee} ),$$  $\PP(U_{\chi})$
 is not contained in any of the proper subspaces
   $$   Fix  (\ga), \ \ for \ \ga  \in G'' \setminus \{ Id\}.$$ 
  Equivalently, $U_{\chi }$ is not contained in a proper eigenspace for the action of an element $g \in G$ on $W_{\chi}^{\vee}$.
  
  We shall say that $V$  is {\bf  projectively general } if every $U_{\chi } \neq 0$ is   projectively general.

\end{defin}

\bigskip

\medskip 
It is convenient to define $E (\chi) : = \sum_{ g \in G} \chi (g^{-1}) g$, and rewrite $ \sA{\chi}$ as a linear  span
 $$ \sA{\chi} = \langle   \sum_{ g \in G} \chi (g^{-1}) \ga g \rangle  _{\ga \in G} = \langle   \sum_{ g \in G} \chi (\ga g^{-1})  g \rangle _{\ga \in G}.$$ 
 
 The subspace $ \sA{\chi} $ maps each element $h \in G$ to the vector 
 $$  (\chi (\ga h^{-1})   ) _{\ga \in G} \in ( \sA{\chi})^{\vee} \subset  \sA^{\vee}.$$
 
 Given elements $h^{-1}, k^{-1}$, they map to the same vector if and only if 
  $$  \chi (\ga h )     =  \chi (\ga k)   , \   \forall \ga \in G \Leftrightarrow  \chi (\ga'   )     =  \chi (\ga'  h^{-1}k)   , \   \forall \ga' \in G .$$
  
  Set $h' : = h^{-1}k$: then our condition amounts to 
   $$  \chi (\ga h' )    = \chi (\ga )   , \   \forall \ga \in G.$$
   
 We specialize our formula to the case where $\ga=1$:   $$  \chi ( h' )     = \chi (1 )  = r : = dim_{\CC} W_{\chi}.$$
 
 \begin{lemma}\label{id}
 If we have a group representation $\rho_{\chi}$ of a finite group $G$, with character $  \chi $,  the condition $  \chi ( h' )  = \chi (1 ) = : r$
 holds if and only $ \rho_{\chi} (h')$ is the identity.
 \end{lemma}
 \Proof
   The eigenvalues of $\rho_{\chi}(h')$,   are 
   $r$ numbers which are roots of unity, hence their real part is at most $1$, and is equal to $1$ if and only if the eigenvalue is $1$.
  Since  the sum of the real parts equals $r$, their real part must equal $1$, therefore    all the eigenvalues are equal to $1$.
  
   Since $h'$ has finite order, we conclude  that  $\rho_{\chi}(h')$ must be equal to the identity.
  
  \qed
   
   We summarize our previous discussion.
   
   \begin{prop}\label{linear}
   a) A  subrepresentation $V \subset \CC[G]$ separates points of $G$  only if, writing $  V =  \sum_{\chi \in Irr}  V_{\chi }, V_{\chi } = W_{\chi } \otimes U_{\chi},$
   
   (i) the representation $$\rho : = \otimes_{ U{\chi} \neq 0} \  \rho_{\chi}$$
   is faithful, i.e. 
  $$  \cap _{ n(\chi) \neq 0} ( Ker \rho_{\chi}) = \{ 1\}. $$
  
 b)  If $V$ is full, or $V$ is linearly general, then  condition (i) is necessary and sufficient.
 
  c) In general, $V \subset \CC[G]$ separates points of $G$ if and only if, setting
  $ \e^*_{\chi} : G \ra  Hom ( U{\chi}, W_{\chi }^{\vee})=  V^{\vee}_{\chi} \subset V^{\vee}$
  $$G_{\chi,V} : =  \{ (g,g') \in G \times G |     \e^*_{\chi} (g) = \e^*_{\chi} (g') \}$$
  $$ \cap_{\chi } G_{\chi,V} = \{ (g,g)\} .$$

       \end{prop}
   
   \Proof
   Since every such representation is contained in a unique minimal full one, in order to prove a) it suffices to prove that (i) 
   is necessary and sufficient  in the case where $V$ is full.
  
   Indeed, we have just shown that two points $h,k$ are separated by $ \sA{\chi}$ if and only if they lie in different cosets of $Ker \rho_{\chi}$.
   
 This shows that  condition (1)  is sufficient, while it is obviously necessary: if $h \in    \cap _{ n(\chi) \neq 0} ( Ker \rho_{\chi})  $
 then $h$ and $1$ are not separated.
 
 The proof of b) in the case where $V$ is linearly general shall  follow from the forthcoming lemma \ref{lem2}.
   
   Finally, c) is a straightforward reformulation of the requirement  that  $\e^*$ separates elements of $G$, since 
  $ \e^* (g) = \e^* (g') \Leftrightarrow   \e^*_{\chi} (g) = \e^*_{\chi} (g') \ \forall \chi$.
  
   \qed
   
   For non full representations, condition (1) is not sufficient, as shown by the following example, which explains in concrete terms
   the notion of a general subrepresentation $V \subset \CC[G]$.
   
   \begin{ex}
Consider the dihedral group $D_n \subset O (2, \RR)$, the group of linear transformations of $\CC$ of the form 
$$ x^j : z \mapsto \zeta^j z,  \ \  x^j y : z \mapsto \zeta^j \bar{z},$$ where $\zeta$ is a primitive  $n$-th root of unity.
We denote by $\rho$ the standard representation $D_n \subset O (2, \CC)$, where, by a suitable change of basis
 $$x (e_1) = \zeta e_1, x (e_2) = \zeta^{n-1} e_2, \ \   y(e_1) = e_2, y (e_2) = e_1.   $$ 
 
 The group algebra is $\sA = \CC[x]/ (x^n-1) \oplus  [\CC[x]/ (x^n-1) ] y$, where $y^2=1, yx= x^{n-1}y$.
 The character $\chi$ of $\rho$  has the property that
 $$ \chi(x^j) = \zeta^j + \zeta^{-j}, \\ \chi (x^j y ) = 0, $$ 
hence 
$$ E (\chi) = \sum_j  ( \zeta^j + \zeta^{-j}) x^j.$$ 

We see easily that $ \sA{\chi} $ has a basis
$$ e_1 : =   \sum_j  \zeta^{-j} x^j, \  e_2 : =   \sum_j  \zeta^{-j} y x^j = y (e_1), $$
$$ e'_2 : = \sum_j  \zeta^j  x^j, \ e'_1 : = \sum_j  \zeta^j  y x^j = y e'_2.$$ 
Moreover $ x(e_1) = \zeta e_1, \  x(e'_1) = \zeta e'_1, \ x(e_2) = \zeta^{-1} e_1, \  x(e'_2) = \zeta^{-1} e'_2.$
 
Hence, for each choice of $(a,b) \in \CC^2 \setminus \{(0,0)\}$ we get a copy of the standard representation inside  $ \sA{\chi} $,
namely, the span $$V = < a e_1 + b e'_1, a e_2 + b e'_2 >.$$

Does $V$ separate points of $D_n$? Calculation shows that 
$$ x^j \mapsto (a \zeta^{-j}, b \zeta^j), \  y x^j \mapsto (b \zeta^j , a \zeta^{-j}).$$

The only possibility for non separation is that $x^j , y x^h$ have the same image, and this holds precisely
when $$ a = b \zeta^{h+j},$$  i.e., exactly when $(a,b)$ is an eigenvector for some reflexion $y x^t$, with eigenvalue $+1$.

 \end{ex}

 \begin{lemma}\label{lem2}
  Let $\chi$ be an irreducible character, and let $V_{\chi } = W_{\chi } \otimes U_{\chi }$, for $U_{\chi }$ a vector subspace of $W_{\chi}^{\vee}$.
  Denote by $G'$ be the factor group $ G / ( Ker (\rho_{\chi})$. $G'$ is separated by $V_{\chi }$ if  $U_{\chi}$ is a linearly general subspace.
 
   \end{lemma} 
 
 \Proof
 $G'$ is separated by $\sA_{\chi} = W_{\chi } \otimes W_{\chi }^{\vee}$, and we have a composition of maps
 $$ G \ra \CC[G]^{\vee} \ra End ( W_{\chi }^{\vee}) \ra Hom (U_{\chi}, W_{\chi }^{\vee}).$$ 
 
 Hence $G'$ is separated by $V_{\chi }$ if and only if there is no pair of distinct elements $\ga, \ga' \in \rho_{\chi} (G) = : G' \subset End ( W_{\chi }^{\vee})$
 such that their restriction to $U_{\chi}$ coincide, equivalently $U_{\chi} \subset Ker (\ga - \ga').$
 
 \qed
   
   \begin{prop}\label{projemb}
   A subrepresentation $V \subset \CC[G]$ yields a map of  $G$ into $\PP (V^{\vee})$  iff $V \neq 0$; it embeds $G$  only if the following properties
   i) and ii) hold.

  i)  Writing $  V =  \sum_{\chi \in Irr}  V_{\chi }, V_{\chi } = W_{\chi } \otimes U_{\chi},$
the representation $$\rho : = \otimes_{ U_{\chi}\neq 0} \  \rho_{\chi}$$
   is faithful.
   
   ii) Let $\chi_1, \dots , \chi_m$ be the irreducible representations for which $ U_{\chi} \neq 0$, and notice  that  $\rho_j$ sends  the centre $Z(G)$ of $G$
to the set of scalar  multiples of the identity,  via a character $\psi_j$ of $ Z(G)$.  Then the characters $\psi_1, \dots, \psi_j$ are affine generators of  
the group $Z(G)^*$ of characters of $Z(G)$
(i.e., $\frac{\psi_j }{\psi_1}$ generate the group $Z(G)^*$).

If moreover $V$ is full, the two conditions i) and ii) are also sufficient.

If $V$ is not full, we obtain a projective embedding of $G$ if $V$ is projectively general (i.e., if each space  $ U_{\chi} \neq 0$  is projectively general).

iii)   In general, $V \subset \CC[G]$ provides a projective embedding of $G$ if and only if, setting
  $ \e^*_{\chi} : G \ra  Hom ( U{\chi}, W_{\chi }^{\vee})=  V^{\vee}_{\chi} \subset V^{\vee}$
  $$G_{\PP(V)} : =  \{ g  \in  G |   \PP (\oplus_{U_{\chi} \neq 0} \e^*_{\chi} )(g) =  Id \  {\rm on  } \ \  \PP (\oplus_{U_{\chi} \neq 0} U_{\chi} )\}$$
consists only of the neutral  element $1 \in G$.
       \end{prop}
 \Proof 
 
 First of all, if $V \neq 0$, there exists $V_{\chi} \neq 0, V_{\chi} = W_{\chi} \otimes U_{\chi}$, and since every element of $G$ yields an 
 element in $GL(W_{\chi})$, $ g (U_{\chi}^{\vee}) \neq 0$.
 
We observe once more that if $V \subset V'$ and $V$ yields a projective  embedding of $G$, then a fortiori also $V'$ does.
Therefore we first deal with the case where $V$ is full, and establish the necessity of conditions i) and ii) in all cases.

Assume that $V$ is full. 
 Now, $h,k$ map to two proportional vectors  inside the dual space $(V )^{\vee}$ if and only if  $1$ and $h' : = h^{-1}k$ map to proportional vectors.
 
 So, wlog, we may assume that $k=1$, and our condition of non injectivity is equivalent to the existence of $ \la \in \CC^*$ such that:
   $$  \chi_j  (\ga h )    = \la  \chi_j (\ga )   , \   \forall \ga \in G, \ \forall j=1, \dots, m.$$

Again here, setting $\ga = 1$, we see that 
 $$  \chi_j  ( h )    = \la  \chi_j (1 ) = \la r_j, \ r_j : = dim (W _{\chi_j}),$$
 while setting $\ga=h$, we get 
  $$  \chi_j  ( h^2 )    = \la  \chi_j (h ) = \la^2   \chi_j (1 )  \Rightarrow   \chi_j  ( h^n )    = \la^n   \chi_j (1 ), \  \forall n \in \NN.$$
  
  If $d$ denotes the order of  $h$, then $\la^d = 1$, so that $\la$ is a root of unity.
At this point the equation $  \chi_j  ( h )    = \la  r_j$ implies that 
$   | \chi_j (h) | = r_j $ , hence all the eigenvalues of $\rho_j (h)$ are equal to $\la$ and $\rho_j (h)$
is a multiple of the identity, $\rho_j (h) = \la Id$.

Hence $\rho_j (h)$  is in the centre of $\rho_j (G)$. Conversely, by Schur's lemma, the centre of $\rho_j (G)$
acts via  scalar multiples of the identity (the scalar being given via a character $\psi_j$ of $Z(\rho_j (G))$.

The same argument in the general case shows iii).

Let us now prove the sufficiency of conditions i), ii).

Since $\rho$ is a faithful representation, if $h$ is not in the centre of $G$, there is a $j$ such that  $\rho_j (h)$  is not in the centre of $\rho_j (G)$. 
Hence $1,h$ are separated.

If instead $h$ is in the centre, $\rho_j (h)$  is  in the centre of $\rho_j (G)$ for all $j$, hence $\rho_j (h) = \psi_j(h) Id$.
By condition ii), if $h \neq 1$, there exist $j,  j'$ such that $ \psi_j(h) \neq \psi_{j'}(h)$, so that $1,h$ are separated.

Conversely, if $\rho$ is not faithful, there exists an element $ h \neq 1$ such that $  \chi_j  ( h )    =   \chi_j (1 ), \ \forall j$,
and $1,h$ are not separated. In instead ii) is not fulfilled, there is an element $h \in Z(G)$ such that $  \chi_j  ( h )    =  \la, \ \forall j$,
hence again $1,h$ are sent to proportional vectors and we do not have projective separation.

{\bf The case where $V$ is not full but projectively general.}

Let $h \neq k$ be elements of $G$. Since $ U_{\chi} \neq 0$  is projectively general, $h,k$ are separated unless their image in $\PP GL (W_{\chi} )$
coincide. If this is the case for all such $\chi$, then $h^{-1} k$ lies in the centre of $G$.

We proceed now as in the full case, using property ii).

\qed
 
 \begin{ex}
Anyone of the full representations of $D_n$ associated to an irreducible 2-dimensional representation separates points of $D_n$.

 The remaining nontrivial  irreducible representations are
just the elements of $D_n^*$: the determinant $\psi_1 : D_n \ra \{\pm 1\}$ is the only one for $n$ odd.
For $n = 2k$ there are also $\psi_2$ and $\psi_1 \otimes \psi_2$, where $\psi_2$ is trivial on reflections, 
and takes values $-1$ on $x$ and on each generator of the subgroup of rotations.

 Therefore, for $n \geq 3$, a full representation $V$  separates points if and only if some irreducible  2-dimensional representation occurs in $V$.
 
 The centre of $D_n$ is trivial for $n$ odd, while for $n$ even it consists of the transformation $  z \mapsto - z$,
 which we shall call $-1$.
 
 Hence for $n$ odd $V$ yields a projective embedding of $D_n$ if some irreducible  2-dimensional representation occurs in $V$.
 If $n$ is even, since in every irreducible  2-dimensional representation $-1$ maps to $- Id$, we need moreover a 1-dimensional 
 character $\psi$ such that $\psi (-1) = 1$.
 
  If $ n = 4 r$, any such $\psi$ does the job, if instead $n=2k$ with $k$ odd, only $\psi_1$
 and the trivial character do the job.
 \end{ex}

 For the reader's convenience, we rewrite the above results in special cases.  
 
\begin{cor}\label{corAb}
Assume that $ V \subset \CC[G]$ is a (left) $\CC[G]$-submodule, $V = V_{\chi_1} \oplus \dots \oplus V_{\chi_m} $. Then

I) If $G$ is Abelian, $V$ separates points of $G$ if and only if the characters $ \chi_1, \dots, \chi_m$
generate the group $G^*: = Hom (G, \CC^*)$ of characters of $G$. 

II) If $G$ is Abelian, $V$ yields and embedding of  $G$ inside $\PP(V^{\vee})$  if and only if the characters $ \chi_1, \dots, \chi_m$
affinely generate the group $G^*$ of characters of $G$. 

III) If $G$ is simple and  $V \neq 0$ is full, then $V$ separates points of $G$ and yields an embedding of  $G$ inside $\PP(V^{\vee})$.
Similarly, if $V \neq 0$ is linearly general, respectively projectively general..
\end{cor}
 \Proof
 In the first two cases we can apply the previous proposition \ref{projemb} since $V$ is full, $Z(G) = G$, and $\rho : =  \chi_1 \times \dots \times  \chi_m$ embeds $G$.
 
   In   case III), $Z(G) = 1$ and every nontrivial homomorphism of $G$ is injective, so i) and ii) hold. 
 
 \qed

 \subsection{Embedding $(G \times G)/G$} 
 
 In general, if we have a homogeneous space $G/H$, which we view as the family of left cosets $gH$ of $H$, and a representation $V \subset \CC[G]$,
 the  map $ \e^*: G \ra V^{\vee} $ does not factor through the quotient map $G \ra G/H$.
 However, if we consider the subspace $V^H$ of $H$-invariants, composing $\e$ with the projection $ V^{\vee} \ra  (V^H) ^{\vee}$
 we obtain a well defined map. 
 It is in general not so easy to see when do we get an embedding of $G/H$, or a projective embedding of $G/H$, 
 except of course in the case where $H$ is a normal subgroup, and $V^H$ can be split as 
 a direct sum of summands $V^H_{\chi}$, where we only take irreducible characters $\chi$ such that $H \subset Ker (\rho_{\chi})$.
 In this way we can apply the result of the previous subsection,
viewing  $V^H$  as a representation  of the quotient group $G' = G/H$.

We shall treat now the following geometrically interesting   case.

 Assume  that we let $G \subset (G \times G)$ be the diagonal subgroup $\{ (g,g) | g \in G\}$, graph of the identity map of $G$.
 
 We consider the homogeneous space $ G^0: = ( G \times G) / G$, the space of left cosets of the diagonal subgroup,
 and we are given two representations $V_1 \subset \CC[G],  V_2 \subset \CC[G]$.
 
 As in the previous subsection $(G \times G)$ maps into  
 $V_1^{\vee} \otimes V_2^{\vee}$ while $G^0$ maps into  the subspace of invariants
 $(V_1^{\vee} \otimes V_2^{\vee})^G = ((V_1 \otimes V_2)^G )^{\vee}$.
 
 We ask once more whether $V_1, V_2$ produce in this way an affine, respectively a projective embedding of $G^0$.
 
 Set for typographical convenience $V_1 = : V(1) , V_2 = : V(2)$.
 
We have first of all by Schur's lemma a decomposition 
$$ (V(1) \otimes V(2))^G  = \oplus_{\chi}  (V(1)_{\chi}  \otimes V(2)_{\chi^{\vee}})^G \subset  \oplus '_{\chi}  (\sA_{\chi}  \otimes \sA_{\chi^{\vee}})^G$$ 
where the symbol $\oplus '_{\chi}$ denotes summing over the (irreducible) characters $\chi$ such that
$V(1)_{\chi}  \otimes V(2)_{\chi^{\vee}} \neq 0$.

 From the abstract point of view, we have  the decomposition
 
 $$ (\sA_{\chi}  \otimes \sA_{\chi^{\vee}})^G  =   [(W_{\chi} \otimes W_{\chi}^{\vee}) \otimes (W_{\chi}^{\vee}. \otimes W_{\chi})]^G = W_{\chi}^{\vee} \otimes W_{\chi},$$
 where we have performed contraction on the first and third terms in view of Schur's lemma.
 
 We derive that   $(V(1)_{\chi}  \otimes V(2)_{\chi^{\vee}})^G \subset (\sA_{\chi}  \otimes \sA_{\chi^{\vee}})^G$
 is then isomorphic to 
 $$ (**) \  U(1)_{\chi}  \otimes U(2)_{\chi^{\vee}} \subset W_{\chi}^{\vee} \otimes W_{\chi}.$$ 
 
  There remains the problem of identifying the right hand side as a space of functions on $G^0 = (G\times G) / G$.
 To be then concrete, let us write as in the previous subsection

$$ (\sA_{\chi}  \otimes \sA_{\chi^{\vee}})  = <  (\sum_{g_1 \in G} \chi (\ga_1 g_1^{-1}) g_1 ) \otimes 
 (\sum_{g_2 \in G} \chi^{\vee} (\ga_2 g_2^{-1}) g_2 ) > _{\ga_1, \ga_2 \in G}$$
 
and observe that the $G$-invariants are just obtained applying the averaging operator $\sum_{\ga} \ga \cdot $, so that 
 
 $$ (\sA_{\chi}  \otimes \sA_{\chi^{\vee}})^G  = <  \sum_{\ga \in G}(\sum_{g_1  \in G} \chi (\ga_1 g_1^{-1}) \ga g_1 ) \otimes 
 (\sum_{g_2 \in G} \chi^{\vee} (\ga_2 g_2^{-1}) \ga g_2 ) > _{\ga_1, \ga_2 \in G} =$$
 $$ = <  (\sum_{g_1, g_2 \in G}  [ \sum_{\ga \in G} \chi (\ga_1 g_1^{-1} \ga)  \chi^{\vee} (\ga_2 g_2^{-1} \ga) ]   g_1  \otimes g_2   ) > _{\ga_1, \ga_2 \in G} .$$ 
 
 In the same way we can explicitly calculate $(V(1)_{\chi}  \otimes V(2)_{\chi^{\vee}})^G$.
  
  We observe now that  for each element of $G^0$ there is a unique representative of the form $(h,1)$,
  and we evaluate the basis tensors on $(h,1)$,
  obtaining by restriction the following elements in the group algebra
   $$   \sum_{h \in G}  F_{\chi, \ga_1, \ga_2} (h) h : =  \sum_{h \in G}  [ \sum_{\ga \in G} \chi (\ga_1 h^{-1} \ga)  \chi^{\vee} (\ga_2 \ga) ]   h       , \ \forall  \ga_1, \ga_2 \in G .$$
   
 We  have
  $$   F_{\chi, \ga_1, \ga_2} (h) =    \sum_{\ga \in G} \chi (\ga_1 h^{-1} \ga)  \chi^{\vee} ( \ga \ga_2) =    \sum_{\ga \in G} \chi (\ga_1 h^{-1}  \ga \ga_2^{-1})  \chi^{\vee} (\ga)=$$
 $$ =  \sum_{\ga \in G} \chi (\ga_2^{-1}  \ga_1 h^{-1} \ga)  \chi^{\vee} (\ga).$$ 
 Hence the above functions depend only upon $ g : =  \ga_2^{-1}  \ga_1$ and we get 
 $$   \sum_{h \in G}  F_{\chi, g} (h) h : =  \sum_{h \in G}  [ \sum_{\ga \in G} \chi (g  h^{-1} \ga)  \chi^{\vee} ( \ga) ]   h       , \ \forall  g  \in G .$$

  Given two elements $h,k$ the two vectors  $ F_{\chi, g} (h)$ and $ F_{\chi, g} (k)$ are equal 
   (resp. linearly dependent)   if and only if,
  defining $$ L_x \chi (\ga) : =  \chi (x \ga )$$ 
  the following Hermitian scalar products  are equal
  $$ < L_{g h^{-1}} \chi, \chi > = < L_{g k^{-1}} \chi, \chi >  \ \forall g \in G,$$
  respectively they are proportional
  $$ \exists \la \ s.t. < L_{g h^{-1}} \chi, \chi > = \la < L_{g k^{-1}} \chi, \chi >  \ \forall g \in G.$$

This holds true for $ h =k$, in general we set $ x : =  k h^{-1}$, and rewrite the previous conditions
as
$$ < L_{g x } \chi, \chi > = \la < L_{g } \chi, \chi >  \ \forall g \in G,$$
$$ \exists \la \ s.t. < L_{g x } \chi, \chi > = \la < L_{g } \chi, \chi >  \ \forall g \in G.$$
The second condition, taking $g = 1$ implies 
$$ < L_{ x } \chi, \chi > = \la,  < L_{ x^2 } \chi, \chi > = \la^2,  \dots \Rightarrow \ \exists m \ s.t. \  \la^m = 1.$$

We can now use the following lemma in order to  infer that if $h,k$ are not separated, then $\rho_{\chi} (h) = \rho_{\chi} (k)$, respectively 
to infer that if $h,k$ are not projectively separated, then $\rho_{\chi} (h^{-1} k) $ lies in the centre of $G' : = \rho_{\chi} (G).$

\begin{lemma}
Let $G$ be a finite group, $\rho$ an irreducible representation, $\chi$ the character of $\rho$. Then 
$$ < L_x \chi , \chi >  = 1 \Leftrightarrow  x \in ker (\rho).$$
If instead  $ < L_x \chi , \chi >  = \la, \ |\la| = 1$, then $\rho (x)$  is in the centre of $G' : = \rho (G)$.
\end{lemma}
\Proof
Without loss of generality we can assume that $\rho$ is faithful, and we show then that in the first case $x=1$.

By the Cauchy-Schwarz inequality, since $ | \chi| = |  L_x \chi | = 1$ (translation invariance of the counting measure),
we obtain $ |< L_x \chi , \chi > | = 1$ iff the two functions are proportional, $ L_x \chi  = \la  \chi  $.
Now $ < L_x \chi , \chi >  = 1$ implies $\la =1$, so that $\chi(x) = 1$,
hence, by lemma \ref{id} $x$ acts as the identity, as we wanted to show.

Assume now that $ < L_x \chi , \chi >  = \la, \ |\la| = 1$: then the same argument shows that  $ L_x \chi  = \la  \chi  $,
hence $\rho (x)$ is a scalar multiple of the identity and lies in the centre of $G'$.

\qed

We  obtain therefore the simplest type of result, which is sufficient to treat the case where the group $G$ is Abelian. We shall pursue elsewhere the analysis of the general case.

\begin{prop}\label{cosetemb}
Assume that $V_1, V_2$ contain respectively  full subrepresentations $V[1], V[2]$
such that,  letting  $\chi_1, \dots , \chi_m$ be the irreducible representations such that   $  V[1] =  \sum_{j=1, \dots, m}  \sA_{\chi_j }$ ,
the following hold:
      
  i)  the  representation $$\rho : = \otimes_{ {j=1, \dots, m}} \  \rho_{\chi_j}$$
  of $G$  is faithful. 
   
   ii) Let  $\psi_j$ be the character of the centre $Z(G)$ of $G$
such that $\rho_j : =  \rho_{\chi_j}$ sends  $ Z(G)$  to the set of scalar  multiples of the identity,  via the character $\psi_j$. 
 Then the characters $\psi_1, \dots, \psi_j$ are affine generators of  $Z(G)$
( $\frac{\psi_j }{\psi_1}$ generate the group $Z(G)$).

iii) $  V[2] =  \sum_{j=1, \dots, m}  \sA_{\chi_j^{\vee}}$.

Then $V_1\otimes V_2$ yields a projective embedding of $G^0 = (G \times G) / G$.
 \end{prop}

 \begin{rem}\label{isotropy}
 A useful variant is the one where $H$ is a normal subgroup of $G$,
  we let $G' : = G/H$, and we consider the orbit space $\bar{G'} : = (G' \times G)/G$,
  which is in an obvious bijection with $G' \times \{1\}$.
  
  We obtain a projective embedding of the orbit space $\bar{G'} $ each time we are given
  representations $V_1, V_2$ of $G'$ satisfying the three properties stated in proposition \ref{cosetemb}.
  
  Just observe that a representation of $G'$ is in a natural way  a representation of $G$!
  \end{rem}

\section{Canonical maps of Galois coverings and of surfaces isogenous to a product}

\subsection{Galois coverings}
 The first situation in which we can apply the results of the previous section is the case where we have a finite flat Galois 
 morphism $f : X \ra Y$ between projective varieties. For simplicity, we make the assumption that both $X,Y$ are smooth,
 and we denote by $G$ the Galois group, so that $ Y = X/G$.
 
 For every locally free sheaf $\om$ of rank $r$ on $X$ such that the action of $G$ on $X$ extends to $\om$, 
 the direct image sheaf $ \sF : = f_* (\om)$ is a  locally free $\hol_Y$-module and a $\hol_Y [G]$-module and accordingly splits as 
 $$ \sF  = \oplus_{\chi \in Irr(G)}  \sF_\chi.$$
 
 If $y \in Y$ is a point such that $f^{-1} (y)$ is in bijection with $G$, then evaluation at the point $y$ (algebraically, this is the operation of tensoring with 
 the  $ \hol_Y$-module $ \hol_Y/ \sM_y$,
 $\sM_y$ being the maximal ideal of the point $y$))
 yields the surjective homomorphism
  $$ \sF  = \oplus_{\chi \in Irr(G)}  \sF_\chi \ra  \CC[G]^r =  \oplus_{\chi \in Irr(G)}  (\sA_\chi)^r, \ r = rk (\om).$$

In the present article we are particularly focused on the case of the canonical sheaf of $X$,
 $\om = \om_X = \hol_X (K_X)$ (here $r=1$).
 
 The $G$-representation we shall consider from now on is
 $$ V : = H^0 ( X, \hol_X (K_X)) = H^0 ( Y, f_* \hol_X (K_X)) = $$
 $$=  \oplus_{\chi \in Irr(G)} H^0 ( Y, \sF_\chi)=  \oplus_{\chi \in Irr(G)} V_\chi. $$
 
  We want to transform basic questions as: is the canonical map of $X$ an embedding, is it birational,
  into necessary or sufficient conditions on the vector bundles $\sF_\chi$, keeping an eye to the examples we shall discuss.
  
  \begin{prop}\label{criterion}
  Let $f : X \ra Y$ be a finite flat Galois morphism between projective varieties, where $X$ is smooth. 
  Let $G$ be the Galois group, and assume that the stabilizers of points of $x$ are normal subgroups $H$
  (this holds if $G$ is Abelian, or if $f$ is unramified).
  
  The canonical system $|K_X|$  is base-point free iff 
  for each point $y \in Y$ the representation $V$ has non trivial image via the composition $ V \ra  \sF \otimes \CC_y  \ra \CC[G/H]$
  (the last homomorphism being associated to the inclusion of the reduced fibre $f^{-1}(y)_{red} \subset f^{-1}(y)$). 
  
  Moreover, in this case  the canonical  morphism  $\Phi$ of $X$ is  injective if: 
  
  (1) for each point $y \in Y$ the representation $V$ embeds the reduced fibre $f^{-1} (y) \cong  G/H$,
  in particular if the associated representation $ V \ra \CC[G/H]$ fulfills properties (i) and (ii) and is projectively general;
  
  (2) $H^0 ( Y, \sF)$ separates pairs of points of $Y$.
  
  The canonical  map $\Phi$ of $X$ is  birational  if (1) holds for some point $y \in Y$,
  and (2) holds for general points $y, y' \in Y$ (we do not need here the condition that $H$ be a normal subgroup).
  
  The canonical map  $\Phi$ of $X$ is an embedding if: 
  
  (1')  for each point $y \in Y$ the representation $V$ embeds the non reduced fibre $f^{-1} (y)$: this holds in particular if
  $\sF$ is  generated by global sections;
  
  (2') $H^0 ( Y, \sF) \ra \sF\otimes \hol_\eta$ is surjective for each length two subscheme $\eta$ of $Y$.

  \end{prop}
  
  \Proof
  
  The first assertion follows from the first assertion of proposition \ref{projemb}.
  
  Injectivity of $\Phi$ means that pairs of  distinct points $x,x'$ have different image.
  
   In the case where $f(x) = f(x')$ 
  (1) guarantees that this is the case,  and for the second assertion we simply use proposition \ref{projemb}
 (recall  that a full subrepresentation is projectively general!).
 
  In the case where $y : = f(x)$ and $y' := f(x')$ are distinct points,
  we use (2) implying that there is a section  $ s \in H^0 ( Y, \sF)$ with $s(y)=0, s(y') \neq 0$. To $s$ corresponds a section
  $\s \in H^0 ( X , \hol_X (K_X))$ with $\s(x)=0, \s(x') \neq 0$.
  
  Birationality of $\Phi$ means that two general points  $x,x'$ are separated, and it suffices to show this for
   $x,x' \in f^{-1} (\sU)$, for $\sU$ a Zariski open set of $Y$. If (1) holds for some point $y$, being an open condition,
   it holds on an a Zariski open set $\sU \subset Y$, and this settles the case where $f(x) = f(x')$, whereas the 
   case where $y : = f(x)$ and $y' := f(x')$ are distinct points is handled as before.
   
   To show embedding, one must moreover show that, for each length two subscheme $ \xi \subset X$,
   $H^0 ( X,  \om_X )  \ra \om_X\otimes \hol_\xi$ is surjective, and it suffices to treat the case where
   $\xi$ is supported at one point $x$, once we have stablished the injectivity of $\Phi$.
   The schematic image of $\xi$ is either a length two subscheme $\eta \subset Y$,
   in which case we use (2'); or $\xi$ is contained in a fibre, and we use (1').
   
  \qed
  
  \begin{example}
  Let $X,Y$ be curves, and assume that $x$ is a ramification point of $f$,
  so that, in local coordinates, $t = f (z) = z^m$.
  We see that in this case the stabilizer $H$ of $x$ is cyclic of order $m$, and 
  that the stalk of $\sF$ at $y = f(x)$ contains $1,z, \dots, z^{m-1}$
  but indeed $1,z$ suffice to embed a tangent vector to $X$ at $x$.
  
  This example shows that the condition that $\sF$ is globally generated is too strong.
  \end{example}
  
  \begin{cor}\label{bir}
 Same notation as in proposition \ref{criterion} and assume that $G$ is Abelian and $Y = Q : = \PP^1 \times \PP^1$.
  
  Then the canonical map is birational if, writing $\sF_\chi = \hol_Q(n_\chi, r_\chi)$, 
  
  a) $H^0 (\sF_\chi) \neq 0 \ ( \Leftrightarrow n_\chi, r_\chi \geq 0) $ holds for a set of characters $\chi$ which
  affinely generate $G^*$;
  
  b) there $\exists \chi , \chi' $ such that  $n_\chi \geq 1,  r_{\chi'} \geq 1.$
  \end{cor}
  \Proof
  a) implies that condition  (1) of proposition \ref{criterion} holds.
  
  If  b) holds for $\chi = \chi'$  (there is a $\chi$ with $n_\chi \geq 1, r_\chi \geq 1$) then  $H^0 (\sF_\chi)$ separates two general points $y,y'$.

  Otherwise, by b), there are $\chi, \chi'$ such that $n_\chi \geq 1, r_\chi=0, n_{\chi'} =0,  r_{\chi'} \geq 1,$
   whence $H^0 (\sF_\chi)$ separates general points $y,y'$ with different first projection to $\PP^1$, and $ H^0 (\sF_{\chi'})$
  separates general points $y,y'$ with different second projection to $\PP^1$.
  
  \qed
  
  \subsection{Canonical maps of SIP 's = surfaces isogenous to a product}
  
  In this section and in the sequel  we shall consider surfaces isogenous to a product of unmixed type
  (see \cite{isogenous} for general properties of SIP ' s and proofs for many statements which we shall make).
  
 We have two curves $C_1, C_2$ of respective genera $g_1, g_2 \geq 2$ and a faithful action of a group $G$ on both curves,
 we denote by $C'_i = C_i / G$ and by $p_i : C_i \ra C'_i$ the quotient morphism.
 
 Clearly $G \times G$ acts on the surface $C_1 \times C_2$, and we consider $G \subset G \times G$
 the diagonal subgroup.  The basic assumption is that $G$ acts freely on the product $C_1 \times C_2$ (this boils down to the fact that,
 denoting by $\Sigma_j$ the set of stabilizers for the action of $G$ on $C_j$, $\Sigma_1 \cap \Sigma_2$
 consists only of $1 \in G$).
 
 The quotient $S : = (C_1 \times C_2) / G$ is then smooth and we have a sequence of morphisms
 $$ \pi: C_1 \times C_2  \ra S ,\  p : S \ra C'_1 \times C'_2 = (C_1 \times C_2)/ (G \times G) = : Y$$
 which are uniquely determined by the requirement that $G$ acts faithfully on both curves.

 We denote by $V_j : = H^0 (C_j, \hol_{C_j} (K_{C_j} ))$, so that 
  $$  H^0 (S,  \hol_{S} (K_{S} )) = (V_1 \otimes V_2)^G = \oplus_{\chi \in Irr(G)} (V_{1, \chi} \otimes V_{2, \chi^{\vee}})^G.$$

 Observe that for $y = (y_1, y_2) \in Y$  the inverse image $p^{-1}(y) $ is the 
 quotient of the  $(G \times G)$-orbit of a point $x = (x_1, x_2) \in  C_1 \times C_2$ by the action of the diagonally embedded $G$.
 
 Since we assume that this action is free, either $x$ is a point where $p \circ \pi$ is unramified,
 and then   $p^{-1}(y) \cong (G \times G) /G$,
 or $x$ is a point where $p \circ \pi$ is not  \'etale, equivalently one of the projections $p_j$ is ramified,
 and then if for instance $j=1$,  $p^{-1}(y) \cong (G/H \times G) /G$,
 where $H$ is the  subgroup of $G$ stabilizing $x_1$.
 
 {\bf Assumption N:} we shall assume, as in proposition \ref{criterion}, that $H$ is a normal subgroup for each point with a nontrivial stabilizer.
 
 The projection $p : S \ra Y$ is Galois if and only if $G$ is Abelian, since $G \subset G \times G$ is a normal subgroup
 iff $G$ is Abelian. We can define as usual
 $$  \sF : = p_* \hol_{S} (K_{S} ),  \e_y :   (V_1 \otimes V_2)^G = H^0 (S,  \hol_{S} (K_{S} )) = H^0(Y, \sF) \ra \sF \otimes \CC_y.$$
 
 In this way we are reduced to the situation considered in proposition \ref{cosetemb} and remark \ref{isotropy}.
 
 Indeed, evaluation on the orbit of $(x_1, x_2)$ yields maps 
 $$  (V_1 \otimes V_2) \ra  V'[1] \otimes V'[2] \subset \CC[G/H_1] \otimes \CC[G/H_2],$$
(here one of the subgroups $H_j$ is trivial, by the assumption that the action of $G$ is free)  inducing 
  $$  (V_1 \otimes V_2)^G  \ra  (V[1] \otimes V[2])^G \subset (\CC[G/H_1] \otimes \CC[G/H_2])^G,$$
  where $V[1], V[2]$ are determined by property iii) of proposition \ref{cosetemb}.
  
  Define moreover $H$ to be trivial  when both $H_1, H_2 $ are trivial, otherwise $H$ is the only nontrivial
  normal subgroup $H_j$. Then $V[1], V[2]$ are representations of $G' : = G/H$. 
  Observe finally that $\CC[G/H] = (\CC[G/H_1] \otimes \CC[G/H_2])^G$ is just the space of functions on the reduced fibre $f^{-1}(y)_{red} \subset f^{-1}(y)$).

  \begin{prop}\label{criterion-isog}
  Let $S$ be a surface isogenous to a product of unmixed type and let $p : S \ra Y = C'_1 \times C'_2 $ be the natural   finite flat  morphism.   Keep assumption N.
    
  Then the  canonical system $|K_S|$  is base-point free iff 
  for each point $y \in Y$ $V : = (V_1 \otimes V_2)^G$ has non trivial image  $(V[1] \otimes V[2])^G$ 
  via the composition $ V \ra  \sF \otimes \CC_y  \ra \CC[G/H]$,   (the last homomorphism being associated to the inclusion of the reduced fibre
  $f^{-1}(y)_{red} \subset f^{-1}(y)$).

  Moreover, in this case  the canonical  morphism  $\Phi$ of $S$ is  injective if: 
  
  (1) for each point $y \in Y$,  $V$ embeds the reduced fibre $f^{-1} (y)_{red} \cong  G/H$, and this holds 
  in particular if $V[1], V[2]$ are full and (i), (ii) of prop. \ref{cosetemb} are fulfilled.
    
  (2) $H^0 ( Y, \sF)$ separates pairs of points of $Y$.
  
  The canonical  map $\Phi$ of $X$ is  birational  if (1) holds for some point $y \in Y$,
  and (2) holds for general points $y, y' \in Y$.
  
  The canonical map  $\Phi$ of $X$ is an embedding if: 
  
  (1')  for each point $y \in Y$ the representation $V$ embeds the non reduced fibre $f^{-1} (y)$: this holds in particular if
  $\sF$ is  generated by global sections;
  
  (2') $H^0 ( Y, \sF) \ra \sF\otimes \hol_\eta$ is surjective for each length two subscheme $\eta$ of $Y$.

  \end{prop}
  \Proof
  The proof is the same as the one of proposition  \ref{criterion},
  using proposition \ref{cosetemb} and remrk \ref{isotropy}.
  
  \qed
  
  \begin{rem}
  Since $\pi$ is unramified, in order to verify that the canonical map is a local embedding,
  it suffices to verify  that it gives a local embedding on $C_1 \times C_2$.
  
In particular, for (1'), it suffices to show that   
  $V[1] \times  V[2]$ embeds every tangent vector which is in the kernel of the projection $p \circ \pi$.
  \end{rem}

\section{Surfaces with record winning high canonical degree}

\subsection{Regular SIP 's with $p_g=4,5,6$}

If we have a SIP, then 
$$K^2_S = 8 (1 - q(S) + p_g(S)) = 8 \frac{1}{|G|}(g_1-1) (g_2-1)$$
and, in order to make $K^2_S$ as high as possible, we require that 
$q(S) = g'_1 + g'_2 = 0$.

This means that $C'_i = \PP^1, i=1,2,$  $Y =  \PP^1 \times \PP^1$,
and then
$$ ( 1 + p_g(S)) |G| =  (g_1-1) (g_2-1).$$

Recall that each covering $p_j : C_j \ra \PP^1$ is  determined by the branch points set  $\sB_j : = \{P_1, \dots, P_{r_j}\}$ and,
once a geometric basis of $\pi_1 (\PP^1 \setminus \sB_j)$ is chosen,  by a spherical system of generators
for the group $G$.
A spherical system of generators for $p_1$ is  the choice of elements $\ga_1, \dots, \ga_{r_1} \in G$
such that $\prod_1^{r_1} \ga_j = 1$. 
Similarly for $p_2$ we get $r_2$ points and  $\ga'_1, \dots, \ga'_{r_2} \in G$
such that $\prod_1^{r_2} \ga'_j = 1$.

One denotes by $m_j$ the order of $\ga_j$, and by $m'_j$ the order of $\ga'_j$.

Then Hurwitz' formula yields

$$ g_1 -1 = \frac{1}{2} |G| [ -2 + \sum_1^{r_1}  ( 1 - \frac{1}{m_j})], $$
$$g_2 -1 = \frac{1}{2} |G| \{ -2 + \sum_1^{r_2}  ( 1 - \frac{1}{m'_j})\}, $$ 
hence 
$$ (*) \  4 ( 1 + p_g(S)) =  |G| [ -2 + \sum_1^{r_1}  ( 1 - \frac{1}{m_j})] \{ -2 + \sum_1^{r_2}  ( 1 - \frac{1}{m'_j})\}.$$

These formulae, even if we restrict to the values $p_g=4,5,6$ lead to a nontrivial determination problem.

Our first search is for uniform  Abelian groups $ G = (\ZZ/m)^k$. 

The $m_j, m'_i$ 's are divisors of the exponent $m$ of the group,
$$m = a_j m_j,  \ m = b_i m'_i.$$

Formula $(*)$ rewrites as an equality among integers:
$$ (**) \  4 ( 1 + p_g(S)) =  m^{k-2} [ -2m  + \sum_1^{r_1}  ( m - a_j)] \{ -2 m + \sum_1^{r_2}  ( m - b_j)\}.$$

And also as 
$$ (***) \  4 ( 1 + p_g(S)) =  m^{k} [ -2  + \sum_1^{r_1}  ( 1 - \frac{1}{m_j})] \{ -2  + \sum_1^{r_2}  ( 1 - \frac{1}{m'_j})\} = : m^k A_1 A_2.$$

{\bf Claim 1:  if  $ k \geq 3$ then $m\leq 3$.}

   Since we have a spherical system of generators, $ r_j \geq k +1$.
 If  $r_1, r_2  \geq 5$, in  formula (***) we have  $A_j \geq 1/2$, hence  
 $$ m^k \leq 16 (1 + p_g),$$
 and for $ k \geq 3$ this implies $m \leq 4$. Actually, $k \leq 6$, and for $k =4 $ $ m \leq 3$, for $k=5,6$, $m=2$.
 
If $r_1=4$, then $k=3$ and all the  $a_j$'s are equal to 1 (three elements are a basis, so also their sum has order $m$)
 hence the term in square brackets is  $[2m - 4]$;
hence if $r_1=4, r_2 \geq 5$ we obtain
$$  m^2 (m-2) \leq  4 (1 + p_g) \Rightarrow  m \leq 3.$$
While, for $r_1=r_2 = 4$ we get  $  m (m-2)^2  =    (1 + p_g) $ and this case is not possible for divisibility reasons.

Assume now $m=4, k=3$: then at least four of the $a_j$'s equal $1$, and the terms $[ ..]$, $\{ ..\}$
are at least $6$, so again this case is not possible, and our claim is established.

\qed

{\bf Case i) $k \geq 3$, $m=2,3$.}

Since $2,3$ are prime numbers, $a_j = b_i = 1, \forall i,j$, and then formula $(**)$ rewrites as 
$$ (****) \  4 ( 1 + p_g(S)) =  m^{k-2} [ -2m  + r_1  ( m - 1)] \{ -2 m + r_2  ( m - 1)\}.$$

Since  $k\geq 3$, then:
$$ p_g = 4 \Rightarrow   m=2, k=3,4$$ 
$$ p_g = 5 \Rightarrow  m=3, k =3, \  {\it or} \ m=2, k=3,4,5$$ 
$$ p_g = 6 \Rightarrow  m=2, k=3,4.$$ 

Recall moreover that $r_j \geq k+1$, because $G$ is not $(k-1)$-generated.
Observe moreover that  if $m=3$ we must  have $r_1 = 4$, $ r_2 = 5$, since $2r_j - 6 \geq 2$.

 Consider now the case $m=2$.

If $m=2$, $k=2$ there are at most $3$ non trivial elements in $G$ hence we cannot get a free action.
Else, we have 
$$ (***) \  4 ( 1 + p_g(S)) =  2^{k-2} [ - 4   + r_1 ] \{ -4 + r_2  \},$$
leading to many possibilities:
$$ p_g = 4 ,  m=2, k=3, r_1 = 6, r_2 = 9, \ {\it or} \ r_1 = 5, r_2 = 14,$$
$$ p_g = 4 ,  m=2, k=4,  r_1 = 5, r_2 = 9,$$
$$ p_g = 5 ,  m=2, k=3,  r_1 = 5, r_2 = 16, \ {\it or} \ r_1 = 6, r_2 = 10,  \ {\it or} \ r_1= 7, r_2 = 8,$$
$$ p_g = 5 ,  m=2, k=4,  r_1 = 5, r_2 = 10, \ {\it or} \ r_1 = 6, r_2 = 7,$$
 $$ p_g = 5 ,  m=2, k=5,  r_1 = 5,  r_2 = 7,$$
$$ p_g = 6 ,  m=2, k=3, r_1 = 6, r_2 = 11, \ {\it or} \ r_1 = 5, r_2 = 18,$$
$$ p_g = 6 ,  m=2, k=4,  r_1 = 5, r_2 = 11.$$

If we want the canonical map of $S$ to be birational, then  at least $k+1$ distinct characters must 
appear in $V$, therefore $ p_g \geq k+1$, and we can exclude the second 
 and the fifth case.

{\bf Case ii) $k=2$, $m \geq 3$ a prime number implies $m=3$ or $m=5, p_g = 4,  6$.}

Here $ \  4 ( 1 + p_g(S)) =   [ ( r_1-2)  ( m - 1) -2] \{ ( r_2-2)  ( m - 1) -2\},$
hence $m \leq 8$. Since $m$ is odd, we can rewrite as 
 $$ 1 + p_g(S)  = [ ( r_1-2)  \frac{ m - 1}{2} -1] \{ ( r_2-2)   \frac{ m - 1}{2} -1 \}.$$

For $m=7$ the terms in brackets are $2$ or $5$ or $8$, or more,  hence 
the equation is not solvable by divisibility conditions. For $m=5$ the terms in brackets are 
$1, 3, 5,  7 \dots$,  so the only solutions are
 $p_g=4, r_1 = 3, r_2 = 5$ and $p_g=6, r_1=3, r_2=6$. For $m=3$ there are  several  solutions $ ( r_1-3)  ( r_2-3) =  1 + p_g$:
either $r_1 = 4$ and $r_2 = 4 + p_g$, or  $p_g=5$ and $r_1 = 5, r_2 = 6$.

\bigskip

\subsection{Examples with $m=2$}
We assume in the following for the sake of simplicity that $m=2$.
Recall moreover that the condition that $G$ acts freely boils down to the condition that the 
set  $\Sigma_1$ of nontrivial elements in the union  of the cyclic subgroups generated by the elements $\ga_j$ does not  intersect the
set  $\Sigma_2$ of nontrivial elements in the union of the cyclic subgroups generated by the elements $\ga'_j$.

\begin{rem}\label{monk=3}
If $k=3$, $m=2$,  then  $G \setminus \{0\}$ has precisely $7$ elements, and $\Sigma_j$ has at least $3$ elements,
hence $|\Sigma_j| = 3,4$ and the cardinalities $|\Sigma_j| $ cannot both equal $4$.

If $|\Sigma_1| = 3$, then the monodromy elements $\ga_j$ are just three basis vectors counted with multiplicity,
and these multiplicities must be even, hence $r_1$ must be even $\geq 6$. Similarly if $|\Sigma_2| = 3$.

Hence, if $|\Sigma_1| = 3$, $|\Sigma_2| = 4$, the first spherical generating system $\Sigma_1$  is of the form, for a suitable basis,
(we use here the classical notation $\ga^r$ to denote that the  vector $\ga$ occurs $r$ times)
$$  e_1 ^{a_1}, \  e_2 ^{a_2}, e_3 ^{a_3}, \ a_j \equiv 0 \ (mod \ 2),$$
while $\Sigma_2$  is of the form (here $e : = e_1 + e_2 + e_3$):
$$ e^b , \ (e_1 + e_2)^{b_3},  \ (e_1 + e_3)^{b_2}, \ (e_2 + e_3)^{b_1}, b \equiv 0 \ (mod \ 2), b_1 \equiv b_2 \equiv b_3 \  \ (mod \ 2).$$

For instance, for $r_2=5$, the only possibility is $e^2,  \ (e_1 + e_2),  \ (e_1 + e_3), \ (e_2 + e_3)$.

If $r_2$ is even, then all the $b_j$'s are even, and $r_2 \geq 8$.

\end{rem}

\begin{rem}\label{restr}
Observe that if $r_j= 5$, and $k =  3$, then  $g_j = 3$ and there are at most $3$ characters occurring in $V_1$.
Hence the canonical map of $S$ cannot be birational since condition (ii) 
of proposition \ref{cosetemb} is certainly violated. We can exclude this case.

Moreover, condition (2) in proposition \ref{criterion-isog} suggests  that the easier cases are those for which either one character contributes to
a subspace of dimension $4 = dim H^0(\hol_Q(1,1))$ (here we use the notation already introduced $Q: = \PP^1 \times\PP^1$) or two contribute to
a subspace of dimension $2 = dim H^0(\hol_Q(1,0))= dim H^0(\hol_Q(0,1))$. 
\end{rem}

We are left with the possibly easy case:
$$ p_g = 6 ,  m=2, k=3, r_1 = 6, r_2 = 11, $$

and with the cases:

$$ p_g = 4 ,  m=2, k=3, r_1 = 6, r_2 = 9,$$
$$ p_g = 5 ,  m=2, k=3,   r_1 = 6, r_2 = 10,  \ {\it or} \ r_1= 7, r_2 = 8,$$
$$ p_g = 5 ,  m=2, k=4,  r_1 = 6, r_2 = 7.$$

For the sake of uniformity of treatment we make the assumption that  $r_1=6, k=3$ (the first three cases).

\begin{lemma}\label{V1}
 To a character $\chi \in ((\ZZ/2)^3)^{\vee}$, we associate the dimension $n_1(\chi) = dim (V_{1,\chi})$.
 
In the case $r_1=6$, $| \Sigma_1| = 3$,
 letting $e_i ^{\vee}, i=1,2,3,$ be the dual basis, and letting $ e ^{\vee} := e_1 ^{\vee} + e_2 ^{\vee} + e_3 ^{\vee}$, we have:
$$n_1 (e_i ^{\vee}) = 0, \    n_1 (e_i ^{\vee} + e_j ^{\vee})  = 1, \ n_1 (e ^{\vee}) =   2     \  \ \ \forall i \neq  j , \ i,j \in \{ 1,2,3 \}.      $$ 

If  more generally   $| \Sigma_1| = 3$, and  the generating system is of the form $  e_1 ^{2 A_1}, \  e_2 ^{ 2 A_2}, e_3 ^{ 2 A_3},$ ($A_j \geq 1$), then $ \ \forall i \neq  j , \ i,j \in \{ 1,2,3 \}, $
$$n_1 (e_i ^{\vee}) =  A_i - 1, \    n_1 (e_i ^{\vee} + e_j ^{\vee})  = A_i + A_j -1, \ n_1 (e ^{\vee}) =   A_1 + A_2+ A_3 -1  .$$
\end{lemma}
\Proof
$n_1(\chi) = dim (V_{1,\chi})$ is the genus of the double cover of $\PP^1$ obtained as $C_1 / Ker (\chi)$.

In the first cases we have a covering branched in two points, in the second ones the covering is branched on four points,
in the third on six points. The corresponding genera are then $0,1,2$.

The proof of the second assertion is similar.

\qed

With the same method we determine $V_2$ using the normal form developed in remark \ref{monk=3}. 

\begin{lemma}\label{V2}
In the case $| \Sigma_2| = 4$,  $\Sigma_2$  is of the form 
$$ e^b , \ (e_1 + e_2)^{b_3},  \ (e_1 + e_3)^{b_2}, \ (e_2 + e_3)^{b_1}, b > 0, b \equiv 0 \ (mod \ 2), b_1 \equiv b_2 \equiv b_3 \  \ (mod \ 2).$$ 

Then if $b = 2B \geq 2, b_j = 1 + 2 B_j,$ then
$$n_2 (e_i ^{\vee}) =  B + B_j + B_k, \    n_2 (e_i ^{\vee} + e_j ^{\vee})  = B_i + B_j, \ n_2 (e ^{\vee}) =   B -1,     \   \{i,j,k \} = \{1,2,3 \}.   $$ 
If instead 
$b = 2 \beta \geq 2, b_j = 2 \beta_j \geq 2,$
$$n_2 (e_i ^{\vee}) =  \beta + \beta_j + \beta_k - 1, \    n_2 (e_i ^{\vee} + e_j ^{\vee})  = \beta_i + \beta_j - 1,  \ n_2 (e ^{\vee}) =   \beta -1     \ , \ \{i,j,k \} = \{1,2,3 \}.   $$
\end{lemma}

\begin{cor}\label{V}
In the case where $| \Sigma_1| = 3$, $| \Sigma_2| = 4$, letting 
$$V = \sum_\chi V_{1,\chi}\otimes V_{2,\chi}, \ 
N (\chi) := dim (V_{1,\chi}\otimes V_{2,\chi})$$ we have, in the case $b = 2B \geq 2, b_j = 1 + 2 B_j,$
$$N (e_i ^{\vee}) =  (A_i-1) (B + B_j + B_k), \    N (e_i ^{\vee} + e_j ^{\vee})  = (A_i + A_j -1) (B_i + B_j), \  \{i,j,k \} = \{1,2,3 \}.    $$
$$ N (e ^{\vee}) =  ( A_1 + A_2+ A_3 -1) ( B-1).  $$ 

Whereas, in the case $b = 2 \beta \geq 2, b_j = 2 \beta_j \geq 2,$ we get
$$N (e_i ^{\vee}) =  (A_i-1) ( \beta + \beta_j + \beta_k - 1), \    N (e_i ^{\vee} + e_j ^{\vee})  = (A_i + A_j -1) (\beta_i + \beta_j - 1), $$
$$ N (e ^{\vee}) =  ( A_1 + A_2+ A_3 -1) ( \beta-1),     \   \ \{i,j,k \} = \{1,2,3 \}. $$ 

\end{cor}

\begin{prop}\label{possible}
In the case 
$$ p_g = 6 ,  m=2, k=3, r_1 = 6, r_2 = 11   $$
if the canonical map of $S$ is birational, then $|\Sigma_2| = 4$, $|\Sigma_1| = 3$, $A_i = 1,  \forall i=1,2,3$, 
$b= 4, b_1 = 3, b_2 =3 , b_3 =1$, and 
$$ N (e ^{\vee}) = 2, \ N (e_1 ^{\vee} + e_2 ^{\vee})= 2, \ N (e_1 ^{\vee} + e_3 ^{\vee})= 1, \ N (e_2 ^{\vee} + e_3 ^{\vee})= 1.$$

In  the case
$$ p_g = 4 ,  m=2, k=3, r_1 = 6, r_2 = 9,$$
the canonical map of $S$ cannot be birational.

In the case
$$ p_g = 5 ,  m=2, k=3,   r_1 = 6, r_2 = 10, $$
if the canonical map of $S$ is birational, then 
we must have
$|\Sigma_2| = 4$, $|\Sigma_1| = 3$, $A_i = 1,  \forall i=1,2,3$, 
$b= 4, b_1= b_2 = b_3 =2$, and in this case
$$ N (e ^{\vee}) = 2, \ N (e_i ^{\vee} + e_j ^{\vee})= 1, \ \forall i\neq j.$$

In the case 
$$ p_g = 5 ,  m=2, k=3,   \ r_1= 7, r_2 = 8,$$
the  canonical map of $S$ cannot be birational.

\end{prop}

\Proof
If $ p_g = 6 ,  m=2, k=3, r_1 = 6, r_2 = 11$ since $r_2$ is odd, $|\Sigma_2| = 4$, hence $|\Sigma_1| = 3$, $A_i = 1,  \forall i=1,2,3$, 
and  we get $b = 2B \geq 2, b_j = 1 + 2 B_j,$ $B + B_1 + B_2 + B_3 = 4$. Since we need four characters to appear, we 
must have $B=2$, $B_1 = B_2 =1$, $B_3=0$.
Then $ N (e ^{\vee}) = 2$, $N (e_1 ^{\vee} + e_2 ^{\vee})= 2$, $N (e_1 ^{\vee} + e_3 ^{\vee})= 1$, $N (e_2 ^{\vee} + e_3 ^{\vee})= 1$.

In the case $ p_g = 4 ,  m=2, k=3, r_1 = 6, r_2 = 9,$ again $r_2$ is odd, $|\Sigma_2| = 4$, hence $|\Sigma_1| = 3$, $A_i = 1,  \forall i=1,2,3$, 
and  we get $b = 2B \geq 2, b_j = 1 + 2 B_j,$ $B + B_1 + B_2 + B_3 = 3$: here at most  three characters can appear, so this case should be excluded.

If $ p_g = 5 ,  m=2, k=3,   r_1 = 6, r_2 = 10, $ again if $|\Sigma_1| = 3$, $|\Sigma_2| = 4$,
$A_i = 1,  \forall i=1,2,3$, then $\beta + \sum_j \beta_j = 5$,
hence there occur four different characters only if $\beta=2, \beta_j = 1 \ \forall j$, and then
$ N (e ^{\vee}) = 2, \ N (e_i ^{\vee} + e_j ^{\vee})= 1, \ \forall i\neq j.$ 

If instead $|\Sigma_1| = 3$, $|\Sigma_2| = 3$, then we use the same formulae as before, except that now one among $ \beta_1, \beta_2,\beta_3$ 
is zero ($\beta \geq 1$ since the other three vectors are linearly dependent) and only three characters can occur.
In fact, if $\beta_3 = 0$, we should have $\beta= \beta_1 = \beta_2 =2$, a contradiction. 

If instead $|\Sigma_1| = 4$, $|\Sigma_2| = 3$, the fact that  $r_1 = 6, r_2 = 10$
 contradicts remark \ref{monk=3}.
 
 In the case 
$ p_g = 5 ,  m=2, k=3,   \ r_1= 7, r_2 = 8,$
let us replace the order and treat $r_1= 8, r_2 = 7.$
Since $r_2$ is odd, $|\Sigma_2| = 4$, hence $|\Sigma_1| = 3$, $A_1 = 2, A_2 = 1, A_3 =1$.
We have $b = 2B \geq 2, b_j = 1 + 2 B_j,$ $ B + B_1 + B_2 + B_3 = 2$.
If  $B=2$,  then $ N (e ^{\vee}) =  3 (B-1) = 3 $: since  $p_g=5$, at most three characters can occur.
If instead $B=1$, then only one $B_j =1$, and there are two cases: $B_1 = 1$, or $B_2 =1$.
In both cases we see that  only three characters occur (with dimensions $ 2,  2 ,1$). 

\qed

\begin{thm}\label{birat}
In the two possible  cases of proposition \ref{possible}

I) $ p_g = 6 ,  G = (\ZZ/2)^3,  r_1 = 6, r_2 = 11   $,
with generating systems $e_1^2, e_2^2, e_3^2$, respectively 
$e^4,  (e_2 + e_3)^3,  (e_1 + e_3)^3, (e_1 + e_2),$ we have
$$ N (e ^{\vee}) = 2\cdot 1, \ N (e_1 ^{\vee} + e_2 ^{\vee})= 1 \cdot 2, \ N (e_1 ^{\vee} + e_3 ^{\vee})= 1 \cdot 1, \ N (e_2 ^{\vee} + e_3 ^{\vee})= 1 \cdot 1,$$ 
and the canonical map of $S$
is a non injective birational morphism  onto a  canonical image $\Sigma$ of degree $56$ in $\PP^5$.

II) $ p_g = 5 ,   G = (\ZZ/2)^3,   r_1 = 6, r_2 = 10, $
with generating systems 
$e_1^2, e_2^2, e_3^2$, respectively 
$e^4, (e_2 + e_3)^2,  (e_1 + e_3)^2, (e_1 + e_2)^2,$ we have 
$$ N (e ^{\vee}) = 2 \cdot 1, \ N (e_i ^{\vee} + e_j ^{\vee})= 1 \cdot 1, \ \forall i\neq j,$$
and the canonical map of $S$
is a non injective birational morphism  onto a  canonical image $\Sigma$ of degree $48$ in $\PP^4$,
for general choice of the branch points for $C_2 \ra \PP^1$.
\end{thm}

\Proof

Let for convenience $u_{ij}=0$ be  the reduced divisor on $C_1$ corresponding to the inverse image of the two points with local monodromy $e_k$,
for $\{1,2,3\} = \{i,j,k\}$.

Similarly define $\s_{ij}=0$ as the reduced divisor on $C_2$ corresponding to the  points with local monodromy $e_i + e_j$,
respectively  $\s=0$ as the reduced divisor on $C_2$ corresponding to the  points with local monodromy $e = e_1 + e_2 + e_3$.
Define $v_{ij}: = \s  \s_{ij} $ and $ v_{123} : = \prod_{i \neq j} \s_{ij} $.

Let us first treat  case I) : $p_g=6$. Here we can write the canonical system as
$$  V = z_{123} H^0(\hol_Q(1,0)) \oplus z_{12} H^0(\hol_Q(0,1)) \oplus \CC z_{13} \oplus \CC z_{23}=$$ 
 $$  V = v_{123} H^0(\hol_Q(1,0)) \oplus u_{12} v_{12} H^0(\hol_Q(0,1)) \oplus \CC u_{13} v_{13} \oplus \CC u_{23} v_{23}.$$

{\bf First step:} The canonical system has no base points.

Indeed,  the base locus is defined by the equations:
$$  v_{123} =   u_{12} v_{12} =   u_{13} v_{13} =  u_{23} v_{23}=0, $$
equivalently 
$$  \prod_{i \neq j} \s_{ij} =   u_{12} \s \s_{12} =   u_{13} \s \s_{13} =  u_{23} \s \s_{23}=0, $$ 
It is crucial to observe that  the divisors of the $\s_{ij} $'s and of $\s$ are pairwise disjoint, 
and the same occurs for the divisors  of the $u_{ij}$'s.

Hence the locus $\s =0$ is away from the base locus, and if $\s_{ij} =0$, then we use that $u_{ik} =0,
u_{jk} =0$ defines the empty set. This proves our first claim. 

Moreover, we see immediately that, since the four above characters are affine generators, condition (1)
of proposition \ref{criterion-isog} is fulfilled for a general point of $Q$; likewise we see that (2) holds 
(see corollary \ref{bir} for a similar argument).

To see that the canonical map is not injective, consider the curve $D = div(\s )$, consisting of the inverse image of $4$ points on the second copy of $\PP^1$.
Here the canonical map sends each component of $D$ to $\PP^1$, hence the canonical map has degree at least $4$ (indeed $16$
since $C_1/e \ra C'_1 = \PP^1 $ has degree $4$).

Let  us now treat the more difficult case II) : $p_g=5$.
In this case a basis of $V$ is provided by 
$$  V = z_{123} H^0(\hol_Q(1,0) \oplus \CC z_{12}  \oplus \CC z_{13} \oplus \CC z_{23}=$$ 
 $$  V = v_{123} H^0(\hol_Q(1,0) \oplus \CC u_{12} v_{12}  \oplus \CC u_{13} v_{13} \oplus \CC u_{23} v_{23}.$$
Here  $u_{ij} $, $v_{ij}= \s  \s_{ij} $ are exactly as before. Hence the same argument shows that the canonical system has no base points.

Likewise, condition (1)
of proposition \ref{criterion-isog} is again fulfilled for a general point of $Q$; to see that (2) holds
we observe that two general points $y,y'$ are separated if they have a different first coordinate $y_1 \neq y'_1$.

Assuming that $y_1 =  y'_1$, then the fibre over a general such point is just $C_2$, and the restriction of the canonical map
$\Phi$ to $C_2$  is projectively equivalent 
to the map $F : C_2 \ra \PP^3$ given by
$$ (v_{123}, v_{12},  v_{13}, v_{23} ) =  (\s_{12}  \s_{13} \s_{23} , \s \s_{12}, \s  \s_{13}, \s \s_{23} ).$$

The map $F$ is equivariant with respect to the action of $G$ on $C_2$ and the action on $\PP^3$
given by the indicated  characters of $G$.

Therefore, if we let $s : \PP^3 \ra \PP^3$ be the map which squares the coordinates,
$s \circ F$ induces a map of $ F' : \PP^1 =C_2 /G \ra \PP^3$, given by 
$$  (\s_{12}^2  \s_{13}^2 \s_{23}^2 , \s^2  \s_{12}^2,  \s^2   \s_{13}^2,  \s^2  \s_{23}^2 )=
 (P_{12}  P_{13} P_{23} , P P_{12}, P  P_{13}, P P_{23} ), $$
 where $P, P_{ij}$ are polynomials on $\PP^1$ (of respective degrees $4,2$) vanishing on the branch points
 which are the  respective images
 of $\s=0 , \s_{ij}=0$.

A projection of $F'$  is the rational map 
$$ F'' =  (P_{12},  P_{13}, P_{23} ),$$
given by three polynomials of degree $2$. 
 Either $F''$ is an embedding of
$\PP^1$ as a conic, and our claim is established, or there is a linear relation among $P_{12},  P_{13}, P_{23}$.
This is however impossible if the three quadratic polynomials   $P_{12},  P_{13}, P_{23}$ are taken to be general.

\qed

\begin{rem}
There remain the cases (with $m=2$)
$$ p_g = 5 ,  m=2, k=4,  r_1 = 6, r_2 = 7,$$
$$ p_g = 6 ,  m=2, k=4,  r_1 = 5, r_2 = 11.$$

In the latter  case (this was the first example found) the canonical map is again a non injective birational morphism for the following choice of generating systems:
$e_1, e_2, e_3, e_4 , e$ and $(e_1 + e_2), (e_1 + e_3)^2, (e_1 + e_4),(e_2 + e_3), (e_2 + e_4), (e_3 + e_4), 
(e_1 + e_2 + e_3), (e_1 + e_ 3 + e_4), (e_1 + e_2+ e_4)^2$.
There occur here $5$ characters, with dimensions $2,1,1,1,1$, so that calculations here are similar to case II) of theorem  
\ref{birat}, but slightly more complicated.
\end{rem}

\section{New double canonical surfaces}

Let us now consider the case $G = (\ZZ/2)^3$, $p_g=4, q=0$, $ r_1 = 6, r_2 = 9$.
Therefore we can use the standard argument yielding 
generating systems 
$e_1^2, e_2^2, e_3^2$, respectively 
$e^{2B}, (e_2 + e_3)^{1 + 2B_1}, (e_1 + e_3)^{1 + 2B_2}, (e_1 + e_2)^{1 + 2B_3}.$
Here, $ B + B_1 + B_2 + B_3 = 3$, hence the only choice in order that the image is not a rational surface
is easily seen to be $B=  B_1 = B_2=1$.

The only characters which occur are the three characters $e_i ^{\vee} + e_j ^{\vee}$, and we get 
$$ \ N (e_1 ^{\vee} + e_2 ^{\vee})= 1 \cdot 2, \ N (e_1 ^{\vee} + e_3 ^{\vee})= 1 \cdot 1, \ N (e_2 ^{\vee} + e_3 ^{\vee})= 1 \cdot 1.$$
It is now immediate to observe that the involution $(e,1) \in G \times G$ ($e =e_1 + e_2 + e_3$) acts on the canonical system as
multiplication by $ + 1$.

Let $Z$ be the quotient surface $S / (e,1)$, so that the canonical map $\Phi$ of $S$ factors through $Z$.
$Z$ is a $ (\ZZ/2)^2$ covering of $Q$ with branch loci of bidegrees $(2,3), (2,3), (2,1)$
made of unions of vertical and horizontal lines, so $Z$ has $14 = 6 + 6 + 2$ nodes,
which are the images of the $14$ isolated fixed points for the involution $(e,1)$ on $S$.
This involution has moreover a fixed divisor, consisting of the reduced inverse image $F'_1 + F'_2$ of the two points 
in the second copy of $\PP^1$ where the local monodromy is $e$.

Since the involution acts as the identity on $H^0( K_S)$, it follows that the isolated fixed points 
occur in the  base locus of the canonical system
with even multiplicity,
while  the  fixed curves occur in the base locus with odd multiplicities $k_1, k_2$.  

Hence if $ q: S \ra Z$ is the quotient morphism,  we claim that $K_S = q^*(K_Z) + F_1 + F_2$.
In fact $40 = K_S^2 = q^*(K_Z + \frac{ k_1}{2} F_1 + \frac{ k_2}{2}F_2)^2$,
with $k_1, k_2 \geq 1$, hence, since
  $K_Z^2=12$, $F_j K_Z =4$ (the two fibrations on $Z$ have respective genera $g_1 = 4,g_2 = 3$,
  $ 20 = 12 + 4 (k_1 + k_2) \Rightarrow k_1 = k_2 =1$.

As in \cite{singular} we see that  
$$  H^0 (Z, \hol_Z (K_Z)) =   H^0 (Q, \hol_Q (0,1)) u_{12} \s_{12} \oplus \CC u_{13} \s_{13}  \oplus \CC u_{23} \s_{23}  ,$$
here again $z_{ij} = u_{ij} \s \s_{ij} $, and the canonical system of $Z$ has no base points.

Moreover, the canonical map of $Z$ is birational if the three quadratic forms $u_{ij}$ are linearly
independent (see proposition 4, pages 104 and 105 of \cite{singular}).

We summarize our discussion in the following

\begin{thm}\label{doublecan}
Surfaces isogenous to a product with $G = (\ZZ/2)^3$, $p_g=4, q=0$, $ r_1 = 6, r_2 = 9$,
and generating systems
$e_1^2, e_2^2, e_3^2$, respectively 
$e^{2}, (e_2 + e_3)^{3}, (e_1 + e_3)^{3}, (e_1 + e_2)^{1}$
yield a connected component of the moduli space of the minimal surfaces 
of general type with $p_g=4, q=0, K^2_S = 40$.

For each of them there is a nontrivial involution $\psi : = (e,1)$ which acts  trivially on the canonical
system, and with quotient a  14-nodal  surface $Z$ with $K_Z^2= 12$, with $K_Z$ ample.

Thus the canonical map of $S$ factors through
 the canonical map of $Z$, which  is a morphism, and is birational onto its
image for general choice of the six branch points of $C_1 \ra  C_1 / G = \PP^1$.

The fixed locus of $\psi$ consists of $14$ isolated points, and two curves which form 
 the reduced inverse image $F'_1 + F'_2$ of the two points 
in the second copy of $\PP^1$ where the local monodromy is $e$.

The canonical system of $S$ has base locus  equal to the divisor $F'_1 + F'_2$.

\end{thm}
\Proof
In \cite{isogenous} it was proven that SIP's with given ramification data form a connected component of the moduli space.

The base points of $|K_S|$ come from the fact that 

$$  H^0 (Z, \hol_Z (K_Z) ) =   H^0 (Q, \hol_Q (0,1))  z_{12} \oplus \CC z_{13} \oplus \CC z_{23} ,$$
and $z_{ij} = u_{ij} \s \s_{ij} $, 
so that we get $\s = 0$ as base locus, and the rest was already shown.

\qed

{\bf Ackowledgement}

The author gratefully acknowledges the support of the 
 ERC Advanced grant n. 340258, `TADMICAMT',
 and of KIAS Seoul, where he has been Research scholar in April 2016.

 Thanks to  Jong Hae Keum and Yongnam Lee for organizing a Workshop `Algebraic Surfaces and Moduli '
 at KAIST Daejeon, March 2016, where   a first example  was presented.

\end{document}